\documentclass[15pt]{article}

\usepackage{amsmath,amssymb, mathtools, amsthm}

\usepackage{amsmath, amsthm, amsfonts, amssymb, mathtools,makeidx}
\usepackage[dvipdfmx]{graphicx,xcolor}

\newcommand{\N}{{\bf N}}

\newcommand{\R}{{\bf R}}

\newcommand{\cF}{\mathcal{F}}
\newcommand{\cL}{\mathcal{L}}

\newcommand{\tF}{\tilde{F}}

\newcommand{\slog}{L}

\newcommand{\ra}{\rightarrow}

\usepackage{mathptmx}
\usepackage{graphicx}
\usepackage{fancybox}
\usepackage{color}

\Large

\def\a{\alpha}

\newtheorem{exam}{Example}[section]
\newtheorem{thm}{Theorem}[section]

\newtheorem{lem}[thm]{Lemma}
\newtheorem{cor}[thm]{Corollary}
\newtheorem{prop}[thm]{Proposition}
\newtheorem{rem}{Remark}[section]
\numberwithin{equation}{section}

\theoremstyle{definition}
     \newtheorem{defn}{Definition}[section]

\newcommand{\dstl}{\displaystyle}

\newcommand{\wt}[1]{\widetilde{#1}}

\newcommand{\wbs}[1]{\,\overline{\!{#1}}}

\newcommand{\wbml}[1]{\;\!\overline{\:\!\!{#1}\;\!\!}}

\newcommand{\ch}[1]{\chi\kern-.05em\lower1ex\hbox{$\scriptstyle{#1}$}}
\newcommand{\chb}[1]
{\wbml{\chi}\kern-.02em\lower1ex\hbox{$\scriptstyle{#1}$}\;\!}

\newcommand{\bs}{\:\!\!\setminus\:\!\!}

\newcommand{\Rn}{{\bf R}\kern-0.08em\lower-0.75ex\hbox{$\scriptscriptstyle n$}}
\newcommand{\Rd}[1]{{\bf R}\kern-0.1em\lower-0.75ex\hbox{$\scriptscriptstyle{#1}$}}
\newcommand{\Sd}{S\kern0.0em\lower-0.75ex\hbox{$\scriptscriptstyle d\;\!\!-\:\!\!1$}}
\newcommand{\Sn}{S\kern0.0em\lower-0.75ex\hbox{$\scriptscriptstyle n\;\!\!-\:\!\!1$}}
\newcommand{\Cn}{{\bf C}\kern-0.01em\lower-0.75ex\hbox{$\scriptscriptstyle n$}}
\newcommand{\el}[1]{{\ell}\kern0.02em\lower-0.75ex\hbox{$\scriptscriptstyle{#1}$}}
\newcommand{\EL}[1]{{L}\kern0.0em\lower-0.75ex\hbox{$\scriptscriptstyle{#1}$}}
\newcommand{\EA}[1]{{A}\kern0.0em\lower-0.75ex\hbox{$\scriptscriptstyle{#1}$}}
\newcommand{\Tst}{T\kern0.0em\lower-0.58ex\hbox{$\scriptscriptstyle *$}}

\newcommand{\wtM}[1]{{\:\wt{\;\!\!\!M\;\!\!}}
\kern0.1em\lower-0.9ex\hbox{$\scriptstyle{#1}$}}
\newcommand{\wbM}[1]{{\;\!\wbs{\:\!\!M\:\!\!}}
\kern0.14em\lower-0.9ex\hbox{$\scriptstyle{#1}$}}

\title{On  the  critical   Caffarelli-Kohn-Nirenberg type inequalities involving  super-logarithms }
\author{ Hiroshi Ando, Toshio Horiuchi, Eiichi Nakai }

\date{}

\begin{document}
\maketitle

\begin{abstract}
We  establish 
the  Caffarelli-Kohn-Nirenberg type inequalities  involving{ super-logarithms (infinitely iterated logarithms).}  
As a result   the critical  Caffarelli-Kohn-Nirenberg type inequalities will be improved, and in certain cases the best constants will be  discovered.
\footnote{2010{\it Mathematics Subject Classification.} Primary 35J70; Secondary 35J60.
\\
 This research was partially supported by Grant-in-Aid for Scientific Research (No. 20K03670, No. 21K03304).}


\end{abstract}

\section{Introduction }
The  main purpose of the present paper is 
to improve the  Caffarelli-Kohn-Nirenberg type inequalities in the critical case studied in \cite{hk3}, 
which are  abbreviated as  the critical CKN-type inequalities.
\par\medskip
We briefly   review the  classical  CKN-type inequalities.
{
The CKN-type inequalities were introduced by \cite{CKN} as  multiplicative interpolation inequalities, but here we refer to the simple weighted Sobolev inequalities. There is a great deal of research in that case alone. In \cite{hk3} we also investigated  the CKN-type inequalities systematically.}
Assume  that $n \ge 1$, $1<p\le q <\infty$ and  $ 0 \le 1/p -1/q \le1/n$.
\par\noindent
$(1)$  In the  non-critical case that $\gamma\neq 0$, the CKN-type  inequalities have the following form:
\begin{equation}\int_{\R^n}|\nabla u(x)|^p |x|^{p(1+\gamma )-n}\, dx\ge
 S^{p,q;\gamma}\bigg(\int_{\R^n}|u(x)|^q |x|^{\gamma q-n}\, dx\bigg)^{p/q}, \; u \in 
C_{\rm c}^{\infty}(\R^n\setminus \{0\}). 
\label{1.4}
\end{equation}
Here $S^{p,q;\gamma}= S^{p,q;\gamma}(\R^n) $ is   the  best constant and given by the following variational problem:
\begin{align}
 S^{\:\!p\;\!\!,q\:\!;\:\!\gamma}
 &=\inf\{{\;\!}E^{\:\!p\;\!\!,q\:\!;\:\!\gamma}[u] : u\in
 C^{\:\!\infty}_{\rm c}(\R^{\;\!\!n}\bs{\;\!\!}\{0\}{\:\!\!}){\;\!\!}\bs{\;\!\!}\{0\}\}, \label{best}
\end{align}
where
\begin{equation}
 E^{\:\!p\;\!\!,q\:\!;\:\!\gamma}[u]
 =\dfrac{\int_{\R^n}|\nabla u(x)|^p |x|^{p(1+\gamma )-n}\, dx}
 {\bigg(\int_{\R^n}|u(x)|^q |x|^{\gamma q-n}\, dx\bigg)^{p/q}}{\qquad}\mbox{ for }u\in
 C^{\:\!\infty}_{\rm c}(\R^{\;\!\!n}\bs{\;\!\!}\{0\}{\:\!\!})\setminus \{0\}.
\end{equation}
 We  also define the  radial best constant  $S^{\:\!p\;\!\!,q\:\!;\:\!\gamma}_{\rm rad}=S^{\:\!p\;\!\!,q\:\!;\:\!\gamma}_{\rm rad}(\R^n)$ as follows.
\begin{equation}
S^{\:\!p\;\!\!,q\:\!;\:\!\gamma}_{\rm rad}
 =\inf\{{\;\!}E^{\:\!p\;\!\!,q\:\!;\:\!\gamma}[u] :  u
 \in C^{\:\!\infty}_{\rm c}(\R^{\;\!\!n}\bs{\;\!\!}\{0\}{\:\!\!})^{}{\!}
 \bs{\;\!\!}\{0\}, \; {  u \text{ is  radial}}{\;\!}\}. \label{symmetric}
\end{equation}
\par\medskip\noindent $(2)$ In the critical case that $\gamma=0$, the CKN-type inequalities have the following form: For $\eta>0$,
let $B_\eta$ be  the  ball $ \{ x\in \R^n : \: |x|< \eta\}$.   
\begin{equation}\label{CKNc}\dstl{
 {\int_{\;\!\!B_{\;\!\!\eta}^{}}\;\!\!}|\nabla u(x)|^{p}|x|^{p-n}{\:\!}dx
 \ge C^{p,q;R}\Big({\int_{\:\!\!B_{\;\!\!\eta}^{}}\;\!\!}
\dfrac{ |u(x)|^{q}}{|x|^n \left(\log(R\eta/|x|)\right)^{1+\:\!q/\;\!\!p'}\;\!\!}{\:\!}dx
 \Big)^{\:\!\!p/\;\!\!q}\,
},\quad  u \in 
C_{\rm c}^{\infty}(B_\eta\setminus \{0\}), \end{equation}
where $p'=p/(p-1)$, $R$ is a positive  number satisfying $ R> 1$  and the best constant $C^{p,q;R}=C^{p,q;R}(B_\eta)$ is given by the  variational problrm:
\begin{align}\label{bestconstantC}
 C^{\:\!p\;\!\!,q\:\!;\:\!R}
 &=\inf\{{\;\!}F^{\:\!p\;\!\!,q\:\!;\:\!R}[u] : u\in
 C^{\:\!\infty}_{\rm c}(B_\eta\bs{\;\!\!}\{0\}{\:\!\!}){\;\!\!}\bs{\;\!\!}\{0\}\}, 
\end{align}
where
\begin{equation}
 F^{\:\!p\;\!\!,q\:\!;\:\!R}[u]
 =\dfrac{\int_{B_\eta}|\nabla u(x)|^p |x|^{p-n}\, dx}
 {\Big({\int_{\:\!\!B_{\;\!\!\eta}^{}}\;\!\!}
{ |u(x)|^{q}}{|x|^{-n }\left(\log(R\eta/|x|)\right)^{-1-\:\!q/\;\!\!p'}\;\!\!}{\:\!}dx
 \Big)^{\:\!\!p/\;\!\!q}}{\,\,\,}\mbox{ for }u\in
 C^{\:\!\infty}_{\rm c}(B_\eta\bs{\;\!\!}\{0\}{\:\!\!})\setminus \{0\}.
\end{equation}
We also define the  radial best constant
\begin{align}
  C^{\:\!p\;\!\!,q\:\!;R}_{\rm rad}
 &
 =\inf\{{\;\!}F^{\:\!p\;\!\!,q\:\!;R}[u] : u
 \in C^{\:\!\infty}_{\rm c}(B_{\eta}^{}{\!}\bs{\;\!\!}\{0\}{\:\!\!})^{}{\!}
 \bs{\;\!\!}\{0\}, \; {  u \text{ is  radial}}{\;\!}\}.\label{symmetric2}
\end{align}
Here we remark  that
   the best constants $C^{\:\!p\;\!\!,q\:\!;\:\!R}$  and 
 $C^{\:\!p\;\!\!,q\:\!;R}_{\rm rad}$ are independent of $\eta>0$. To see this it  suffices  to employ a change  of  variables given by $x=\eta y$.
 Roughly speaking 
if $p\ge n$, the imbedding inequalities (\ref{CKNc}) follow from  a generalized rearrangement argument.
On the other hand if $1 < p < n$, then these are established  by using the
so-called nonlinear potential theory. 
\par 
Recently  M. Sano and F. Takahashi established in  \cite{ST} 
the CKN-type inequality with logarithmic weights in two dimensions using  the nonlinear potential theory. 
In  \cite{ho3,ho4,ho5} we also  revisited the CKN inequality  and extended them to the case with non-doubling weights. 
Interestingly, the  framework in \cite{ho4} not only unifies non-critical and critical cases, but also makes it possible to study 
the critical case in greater depth. For more information on the CKN-type inequalities and relating inequalities see \cite{CKN, m, ho1, ho2}.
On the other hand, 
the  super-logarithms have been introduced and developed by  \cite{AHN-2012,FO,AHN-2014, AHN-2014-2}.
Intuitively, the super-logarithms are defined as  concave functions that grow more slowly than all poly-logarithms and they are realized
as arbitrarily and infinitely iterated logarithm functions. They have  been  applied to  the refinement of  the Hardy inequalities by finding infinitely many missing terms in \cite{AHN-2014, AHN-2014-2} and  the analysis of weighted $L^p$-boundedness of convolution type integral operators in \cite{FO}. 

\par\medskip 
In this  paper we will {improve  the critical  CKN-type inequalities represented by (\ref{CKNc}) by using  not only poly-logarithms  but also super-logarithms as weights instead of $|x|^{p-n}$ on the left side,}
 and in  certain cases  the best constants of them will be determined. 
 To do this we use
 recent results on the  CKN-type inequalities from \cite{ho4} instead of the  nonlinear potential theory,
 and a generalized  rearrangement of  function from \cite{hk3}. 
 
\par\medskip
This paper is  organized in the following way:
\S 2
 provide the basic and essential knowledge of this paper. In the first \S2.1,  super-logarithms are reviewed  from \cite{AHN-2012}.
Next, \S2.2 reviews  CKN-type inequalities  in the new framework  from \cite{ho4}.
Finally, in \S2.3 we  recall a generalized rearrangement of  function  from \cite{hk3}.
After these preparations, { as preliminary results,} \S 3 gives critical CKN-type  inequalities { involving poly-logarithms as weights.
Then,} in \S4, we state the critical CKN-type inequalities involving super-logarithms and also give  proofs of  them.
{Moreover,} \S 5 gives further results.

\section{Preliminaries}

\subsection{{A review on super-logarithms (infinitely iterated logarithms)}}
{
 We recall  poly-logarithm  and  poly-exponential functions as  follows.
\begin{equation}\label{poly-logarithm}\begin{cases}&
\log^0 r=r,  \;  \log^n r=\log(\log^{n-1} r), \qquad \quad  n\in\N,
\\
&
\exp^0(r)=r, \;  \exp^n(r)=\exp(\exp^{n-1}(r)), \;\;  n\in\N.
\end{cases}
\end{equation}
}

First we define two sets of functions{, following \cite{AHN-2012}.}

\begin{defn}\label{defn:cL}
Let $\cL$ be the set of all continuous, increasing and bijective functions $f$
from $(0,\infty)$ to $(-\infty,\infty)$ satisfying 
\begin{equation*}
 \lim_{r\to +0}f(r)=-\infty,
 \quad
 f(1)=0,
 \quad
 \lim_{r\to\infty}f(r)=\infty.
\end{equation*}
\end{defn}
For example, the logarithmic function $\log r$ is in $\cL$.

\begin{defn}\label{defn:cF}
For $a>1$, 
let $\cF_a$ be the set of all continuous, increasing and bijective 
functions from $[a,\infty)$ to itself. 
\end{defn}

If $f\in\cF_a$, then $f(a)=a$ and $\displaystyle\lim_{u\to\infty}f(u)=\infty$.
For a function $f\in\cF_a$, 
let $f^{0}(u)=u$ and $f^{k}(u)=f(f^{k-1}(u))$, $k\in\N$.
Then $f^k$ is also in $\cF_a$.

We define a function $F\in\cF_a$ as
\begin{equation}\label{F}
 F(u)=F_a(u)=a-\log a+\log u \quad (u\ge a).
\end{equation}
Then the relation
\begin{equation}\label{F'}
 (F^{k}(u))'=\frac1{F^{k-1}(u)\cdots F^{1}(u)F^{0}(u)}
\end{equation}
holds.
That is,
\begin{equation}\label{F-int}
 F^{k}(u)=a+\int_a^{u}\frac{dt}{F^{k-1}(t)\cdots F^{1}(t)F^{0}(t)}.
\end{equation}
\begin{defn}\label{defn:sphi}
For $a>1$, let
\begin{equation}\label{tF}
 \tF(u)=\tF_a(u)=a\prod_{k=0}^{\infty}\frac{F^{k}(u)}{a}= {u\prod_{k=1}^{\infty}\frac{F^{k}(u)}{a}} \quad (u\ge a).
\end{equation}
\end{defn}

\begin{defn}[Super-Logarithm]\label{defn:slog}
For $a>1$, let
\begin{equation}\label{slog}
 \slog(r)
 =
 \slog_a(r)
 =
 \int_a^{ar}\frac{1}{\tF(t)}\,dt
 \quad (r\ge 1),
\end{equation}
and let
\begin{equation}\label{slog0}
 \slog(r)
 =-\slog(1/r)
 =
 -\int_a^{a/r}\frac{1}{\tF(t)}\,dt
 \quad (0<r<1),
\end{equation}
where $\tF$ is as in \eqref{tF}.
\end{defn}

Then we have the following.
\begin{thm}[\cite{AHN-2012}, Theorem 2.1]\label{thm:F}
Let $a>1$.
\begin{enumerate}
\item 
The function $\tF$ is in $\cF_a$,
infinitely differentiable
and has the following expression:
\begin{equation}\label{V}
 \tF(u)=\exp(V(u)),
 \quad
 V(u)
 =\log a+\int_a^u
  \left(\sum_{k=0}^{\infty}\frac1{\prod_{j=0}^{k}F^{j}(t)}\right)\,dt.
\end{equation}
Further,
$\left(\dfrac{d}{du}\right)^{k}\dfrac{\tF'}{\tF}$ is bounded 
for each $k\in\{0\}\cup\N$.
\item
The function $\slog$ is in $\cL$, 
differentiable on $(0,\infty)$,
infinitely differentiable except at~$1$,
and, concave on $[1,\infty)$,
Moreover, if $a\ge2$, then $\slog$ is concave on $(0,\infty)$.
\item
For each $n\in\N$,
\begin{equation}\label{SL/log-n}
 \lim_{r\to +0}\frac{\slog(r)}{\log^n(1/r)}=
 \lim_{r\to\infty}\frac{\slog(r)}{\log^n r}=0.
\end{equation}
\item
For $r\ge\exp(a)$,
\begin{equation*}
 \slog(r)
 \le
 \slog(\exp(r))
 \le
 (1+a)\slog(r).
\end{equation*}
\end{enumerate}
\end{thm}
%
\par\medskip

\begin{rem}
We note that the following relations 
\begin{equation}
\lim_{r\ra\infty}\frac{F^k(ar)}{\log^k r} = 1, 
\quad 
  \lim_{r\ra\infty} \frac{F^k(ar)}{\displaystyle\prod_{j=1}^{\infty}\frac{F^{j}(ar)}{a}} = 0, 
\quad 
 \lim_{r\ra\infty} \frac{a+L(r)}{F^k(ar)} = 0, 
\quad 
  \lim_{r\ra\infty}\frac{F^k(a+L(r))}{\log^k (a+L(r))} = 1 \label{relations}
\end{equation}
and
\begin{equation}
\begin{cases}
\dstl  L'(r)= \frac{a}{ \tF(ar)}=\frac{1}{{r}\displaystyle\prod_{j=1}^{\infty}\frac{F^{j}(ar)}{a}} \qquad \quad &(r>1),\\[6ex]
 \dstl (\tF(u))'= \tF(u) \left(\dstl\sum_{k=0}^{\infty}\frac1{\prod_{j=0}^{k}F^{j}(u)}\right) & (u>a).
\\
\end{cases}\label{2.34}
\end{equation}
hold  for any $ k \in \N $ and  $a>1$.

\end{rem}

\subsection{A review on the  general CKN-type inequalities with non-doubling weights}
{$\R_+= (0,\infty)$.}
By $C^{0,1}(\R_+)$ we  denote  the space  of  all Lipschitz continuous functions on $\R_+$.
First we  define a class of weight functions $W(\R_+)$ which  is a slight modification of the one adopted in \cite{ho4} 
to suit our purpose.
\begin{defn}\label{D0}
 Let us  set 
  \begin{equation}W(\R_+)=
\{ w\in C^{0,1}({\R}_+): w>0, \lim_{t\to+0}w(t)=a\, \text{  for  some }\,a\in [0,\infty] \}.\end{equation}

\end{defn}
In  the  next  we define two subclasses of  this   large class.

\begin{defn}\label{D1} Let  us  set 
 \begin{align} &P(\R_+)= \{ w\in W(\R_+) : \,  w^{-1} \notin L^1((0,\eta)) \, \text{ for some } \, \eta >0\}.\\
&
 Q(\R_+) =\{ w\in W(\R_+) :  \, w^{-1}\in L^1((0,\eta)) \, \text{ for any } \, \eta >0 \}.\end{align}
\end{defn}

\begin{exam} 

$t(\log (1+1/t))^\a\in P(\R_+)$ if $\a \le 1$  and 
$t(\log (1+1/t))^\a\in Q(\R_+)$ if $\a > 1$.
\end{exam}

\begin{rem} 
 From  {Definitions \ref{D0} and  \ref{D1}} it  follows   that $W(\R_+)= P(\R_+)\cup Q(\R_+)$ and  $P(\R_+)\cap Q(\R_+)=\emptyset.$
\end{rem}
\par\medskip
We define   functions $f_{\eta}(t)$  and $G_\eta(t)$ on $(0,\eta]$ in order to introduce  variants of the Hardy potentials.
\par\medskip
\begin{defn} \label{Definition2.3}
For $\mu>0, \eta>0$ and  $t\in (0,\eta]$ we set the followings:
\begin{enumerate}\item When $w\in P(\R_+)$,
\begin{equation}f_\eta(t; w,\mu)= \mu+ \int_t^\eta \frac{1}{w(s)}\,ds.\label{fetaP}
\end{equation}
\item
When $w\in Q(\R_+)$,
\begin{equation}f_\eta(t;w)= \int_0^t \frac{1}{w(s)}\,ds.
\label{fetaQ}
\end{equation}
\item
For $w\in W(\R_+)$
\begin{equation}
 G_\eta(t; w,\mu)= \mu+  \int_t^\eta  \frac1{w(s)\; f_\eta(s)}\, {ds}. \label{Geta}
 \end{equation}
 \item $f_{\eta}(t;w,\mu)$ and $f_{\eta}(t;w)$ are abbreviated as $f_{\eta}(t)$.
$ G_{\eta}(t; w,\mu)$ is abbreviated as  $ G_{\eta}(t)$.
\end{enumerate}

\end{defn}
\par\medskip

{
Now  we introduce {\bf the non-degenerate condition (NDC)} as  follows. 
This   is an explicit rewrite of  the condition  which was first  introduced in \cite{ho4} as Definition 3.3 :  }
\begin{defn}[{ the non-degenerate condition}]  Let $\eta>0$  and $w\in W(\R_+)$.  A weight function $w$ is said to satisfy   the non-degenerate condition  if  
\begin{equation}
C_0:=\dstl\inf_{0<t\le \eta} \frac{w(t)}{t}\, f_\eta(t)>0.
\label{H-condition} \end{equation}

\end{defn}
In \cite{ho4} we  have established the $n$-dimensional CKN-type inequalities for  $p>1$.   
 We  cite the following  result as the basic knowledge of this paper.
 Let  $B_\eta$ be  the ball $  \{x\in \R^n:  |x| <\eta\}$. 
\begin{thm}[\cite{ho4}, Theorem 3.1]\label{th4.2} Let $n\in \N$, $1<p\le q <\infty$, $\mu>0$, $\eta>0$ and $ 0 \le 1/p -1/q \le1/n$. 
Assume that $w\in W(\mathbf R_+)$.
Moreover assume  that if $n>1$ and  $p<q$,  $w$ satisfies the non-degenerate condition (\ref{H-condition}).
Then,  we have  the followings:
\begin{enumerate}
\item
There exists a positive  number $C_n=C_n(p,q,\eta,\mu,w)$ such  that
 for any $u\in C_c^\infty(B_\eta\setminus \{0\})$  we  have
\begin{equation}
\int_{B_\eta} |\nabla u(x)|^p w(|x|)^{p-1}|x|^{1-n}\,dx\ge  C_n \left(  \int_{B_\eta} 
\frac{ |u(x)|^q |x|^{1-n}\,dx}{ w(|x|) f_\eta (|x|)^{1+q/{p'}}} \right)^{p/q}\label{6.16'}
\end{equation}
where $p'=p/(p-1)$  and  $f_\eta(t)$ is  given by Definition \ref{Definition2.3}.  
\item 
If $C_n$ is the best constant, then $C_n$ satisfies the followings:
\begin{enumerate}\item
If $n=1$, then $C_1=   S^{p,q;1/p'}(\R)= 2^{p/q-1} S_{ \rm rad}^{p,q;1/p'}(\R)$.
\item 
If  $n>1$, then 
$C_n\ge \min(C_0^p,1)S^{p,q;1/p'}(\R^n)\:$ $(p<q)$; $\; C_n=(1/p')^p\:$ $(p=q)$. \par\noindent 
\item  If $n>1$,  $C_0\ge 1$ and  $1/p'\le \gamma_{p,q}$, then 
\begin{equation} C_n=  S^{p,q;1/p'}(\R^n)=S_{\rm rad}^{p,q;1/p'}(\R^n),\label{2.13}\end{equation}
\end{enumerate}
where $C_0$ is  given by (\ref{H-condition}) and  $\gamma_{p,q}$ is given by 
\begin{equation}\dstl{
 \gamma_{p\;\!\!,q}=\dfrac{n-{\;\!\!}1}{1{\;\!\!}+q/p'\;\!\!}.}\label{gamma}
\end{equation}
Here  $S^{p,q;1/p'}(\R^n)$  and $S_{\rm rad}^{p,q;1/p'}(\R^n)$ are 
{independent of $\eta$. They are called}
the  best constants of  the classical CKN-type inequalities and  
defined in Introduction.
\end{enumerate}
\end{thm}
\par\medskip
\begin{rem}\label{remark2.2} \begin{enumerate}\item
Theorem \ref{th4.2} was established in \cite{ho4} as Theorem 3.1 with $C^1$-weight functions. 
It is an easy task to check  that it remains valid for weight functions in {$W(\R_+)\; (\subset C^{0,1}(\R_+))$} defined in {Definition \ref{D0}.}
\item

If we consider the inequality (\ref{6.16'}) in  the radially symmetric space 
\begin{equation}
C^{\:\!\infty}_{\rm c}(B_\eta\bs{\;\!\!}\{0\}{\:\!\!})_{\rm rad}=
\{ u\in C^{\:\!\infty}_{\rm c}(B_\eta\bs{\;\!\!}\{0\}{\:\!\!}) : u \text{ is radial}\}
\label{radiallysymmetric}
\end{equation}
in place  of $C_c^\infty(B_\eta\setminus\{0\})$, then  the assertion (i) {in Theorem \ref{th4.2}}  holds without the non-degenerate  condition (\ref{H-condition}).
By $C_{n,{\rm rad}}$ we  denote the best constant of the inequality (\ref{6.16'}) in  
$C^{\:\!\infty}_{\rm c}(B_\eta\bs{\;\!\!}\{0\}{\:\!\!})_{\rm rad}$. 
Then the best constant $C_{n,{\rm rad}}$  always satisfies 
$C_{n,{\rm rad}}= S_{\rm rad}^{p,q;1/p'}(\R^n)$  if $n\ge 2$.
In fact, 
by employing  a polar coordinate system, the proof  is done in a quite  similar way  as in  \cite{ho4}.

\end{enumerate}
\end{rem}
{The following lemma gives a sufficient condition for  $1/p'\le \gamma_{p,q}$ in Theorem \ref{th4.2} (ii) (c). }

\begin{lem} \label{lem2.1}
Assume that $n>1$, $\dstl 0\le 1/p - 1/q\le 1/n$ and $1<p\le \dstl{(n+1)}/{2}$. Then, it holds  that
\begin{equation}1/{p'}\le \gamma_{p,q}.
\end{equation}

\end{lem}

\noindent{\bf Proof: } Let $ 1/p - 1/q=s \;( 0\le s\le 1/n)$. Then $ q= p/(1-sp)$.
By (\ref{gamma}) it  suffices to  show that
$$1-1/p\le \frac{ (n-1)(1-sp)}{p(1-s)}.$$
If $n>2$, then   we have 
$$p\le \frac{n-s}{(n-2)s +1}= \frac{1}{n-2}\left(\frac{n+1/(n-2)}{ s+1/(n-2)}-1\right).
$$
The right-hand side takes its minimum  ${(n+1)}/{2}$ at $s=1/n$, hence  the assertion follows.
If $n=2$, then the assertion is also clear.
\qed
\par\medskip

\subsection{ A generalized rearrangement of functions}
We recall  a  rearrangement of functions with respect to general weight functions  instead of Lebesgue measure, which was used in \cite{hk3} to 
establish the validities of ${\;\!}\dstl{
 S^{\:\!p\;\!\!,q\:\!;\:\!\gamma}
}=\dstl{
 S^{\:\!p\;\!\!,q\:\!;\:\!\gamma}_{\rm rad}
}$ 
  and $\dstl{
 C^{\:\!p\;\!\!,q\:\!;R}
}=\dstl{
 C^{\:\!p\;\!\!,q\:\!;R}_{\rm rad}
}{\;\!}$   
under additional conditions. 

\begin{defn}\label{df4.2}\begin{enumerate}
\item For 
$g\in\dstl{
 L^{1}_{\rm loc}(\R^{\;\!\!n}\;\!\!)
}$  and   ${\,}g\ge 0{\,\,}$ a.e. on $\dstl{
 \R^{\;\!\!n}
}$, let us set for  a (Lebesgue) measurable  set  $A$ 
\begin{equation}\dstl{
 \mu_{\:\!\!g}^{}(A)={\int_{A}}d\mu_{\:\!\!g}^{}={\int_{A}}g(x){\:\!}dx
}.\end{equation}
Then 
 $\mu_{\:\!\!g}^{}$ is said to be the measure  determined by $g$. 
\item 
$g$ is said to be  admissible, if  and only if  ${\;\!}g\in\dstl{
 L^{1}_{\rm loc}(\R^{\;\!\!n}\;\!\!)
 \cap C(\R^{\;\!\!n}\bs{\;\!\!}\{0\}{\:\!\!})_{\rm rad}^{}
}$, ${\,}g\ge 0{\,\,}$ on $\dstl{
 \R^{\;\!\!n}\bs{\;\!\!}\{0\}
}{\,}$ and  $g$ is non-increasing with respect to  ${\:\!}r={|x|}{\:\!}$.
{
Here $C(\R^{\;\!\!n}\bs{\;\!\!}\{0\}{\:\!\!})_{\rm rad}^{}=\{ u\in C(\R^{\;\!\!n}\bs{\;\!\!}\{0\}{\:\!\!}) : u \text{ is radial}\}.
$
}
\item 
    For an admissible $g$  and a Borel set $A\subset  {\bf R}^n$ satisfying $ 0<\mu_1(A)<+\infty$, let us define $r_g[A]>0$ by $\mu_g(A)= \mu_g( B_{{r_g}[A]})$. 
   {Here 
    the  ball $B_{{r_g}[A]}= \{x\in \R^{\;\!\!n}: |x|< {r_g}[A]\}$ is said to be the  rearrangement  set of $A$ by $g$.}
\item
For an admissible {$g$} and  ${\,}u:\dstl{
 \R^{\;\!\!n}
}\to\R{\:\!}$,  
 we  set
\begin{align} &\dstl{
 \mu_{\:\!\!g}^{}[u](t)=\mu_{\:\!\!g}^{}(\;\!\!\{{\:\!}|u|>t{\:\!}\}\;\!\!)
 ={\int_{\{|u|\;\!>\;\!t\}}}g(x){\:\!}dx \hspace{0.76cm} \mbox{ for }\; t\ge 0,
}\\ 
&\dstl{
 {\cal R}_{\:\!\!g}^{}[u](x)= {\cal R}_{\:\!\!g}^{}[u](|x|)=\sup\{{\;\!}t\ge 0\,:\,
 \mu_{\:\!\!g}^{}[u](t)>\mu_{\:\!\!g}^{}(B_{|\;\!\!x\;\!\!|}^{}\;\!\!)\}{\,\,\,}
 \mbox{ for }\; x\in\R^{\;\!\!n}\bs{\;\!\!}\{0\}
}.\end{align}
Then  ${\mu_{\:\!\!g}^{}[u]}$  and  ${{\cal R}_{\:\!\!g}^{}[u]}$ are  said to be 
 the  distribution function of $u$ and  the  rearrangement  function  of  $u$ with respect to $g$, respectively. 
\end{enumerate}
\end{defn}

\medskip
Direct from this definition we see the next proposition.

\begin{prop}\label{prop4.3}
Let 
$1\le p<\infty{\:\!}$  and assume that $g$ is admissible. Then,  for ${\,}u:\dstl{
 \R^{\;\!\!n}
}\to\R{\:\!}$, we have  the followings: 
\begin{enumerate}\item

$\dstl{
 \mu_{\:\!\!g^{}}[u](t)=\mu_{\:\!\!g}^{}[{\;\!}{\cal R}_{\:\!\!g}^{}[u]](t){\,\,\,}
 \mbox{ for }t\ge 0.
}$
\item
$\dstl{
 {\cal R}_{\:\!\!g}^{}[|u|^{\:\!p}](x)={\cal R}_{\:\!\!g}^{}[u](x)^{p}{\,\,\,}
 \mbox{ for }x\in\R^{\;\!\!n}\bs{\;\!\!}\{0\}.
}$
\item If 
$u$ is  radially symmetric and  non-increasing with respect to  ${\:\!}r={|x|}{\:\!}$,  then
$$\dstl{
 {\cal R}_{\:\!\!g}^{}[u](x)=u(x) \hspace{1.16cm}
 \mbox{ for \;  a.e. }x\in\R^{\;\!\!n}\bs{\;\!\!}\{0\}.
}$$
\end{enumerate}
\end{prop}
\bigskip
Further we have 
\begin{prop}\label{prop4.4}
Let 
$1\le p<\infty{\:\!}$ and  assume that $g$ is admissible. Then,  for ${\,}u,v:\dstl{
 \R^{\;\!\!n}
}\to\R{\:\!}$, we have  the followings: 
\begin{enumerate}
\item
$\dstl{
 {\int_{\R^n}}|u(x)|^{p}g(x){\:\!}dx={\int_{\R^n}}{\cal R}_{\:\!\!g}^{}[u](x)^{p}g(x){\:\!}dx.
}$
\item
$\dstl{
 {\int_{\R^n}}|u(x){\:\!}v(x)|{\:\!}g(x){\:\!}dx
 \le {\int_{\R^n}}{\cal R}_{\:\!\!g}^{}[u](x){\;\!}{\cal R}_{\:\!\!g}^{}[v](x){\:\!}g(x){\:\!}dx.
}$
\end{enumerate}
\end{prop}
\par 
\medskip
By $u\in C^{0,1}_{\rm c}(\R^n)$ we  denote   the space  of  all Lipschitz continuous functions on $\R^n$ having compact supports.
If $u\in C^{0,1}_{\rm c}(\R^n)$, then $u$ is differentiable for  a.e. $x\in \R^n$ and $|\nabla u|\in L^\infty(\R^n)$.
For an admissible $g$,  we  see that ${\cal R}_{\:\!\!g}^{}[u]$ for $u\in  C_{\rm c}^{0,1}({\R^n})$  is  {radially symmetric and non-increasing}, and  hence ${\cal R}_{\:\!\!g}^{}[u]$ 
is differentiable  for a.e. $x\in \R^n$. Moreover  we have the following. For the proof  see \cite{hk3}, Proposition 4.5.
\par\medskip
\begin{prop}\label{prop4.5}
Let
$1\le p<\infty{\:\!}$ and  assume  that  $g$ is admissible. Then,  for ${\:\!}u\in\dstl{
 C^{0,1}_{\rm c}(\R^{\;\!\!n}\;\!\!)
}{\:\!}$ we have 

$${\int_{\R^n}}|\nabla u(x)|^{p}\dfrac{1}{g(x)^{p\:\!-1}\;\!\!}dx\ge 
\dstl{
 {\int_{\R^n}}|\nabla[{\;\!}{\cal R}_{\:\!\!g}^{}[u]](x)|^{p}\dfrac{1}{g(x)^{p\:\!-1}\;\!\!}dx.
}$$
\end{prop}
From Proposition \ref{prop4.4} and Proposition \ref{prop4.5} we  have  the following:
\begin{prop}\label{prop4.6}
Let
$1\le p<\infty{\:\!}$ and  assume  that  $g$ is admissible. 
Moreover assume  that ${\,}v:\dstl{
 \R^{\;\!\!n}}\to\R_+{\:\!}$ and $ v$ is radially symmetric  and {non-increasing}.
Then,  for ${\:\!}u\in\dstl{
 C^{0,1}_{\rm c}(\R^{\;\!\!n}\;\!\!)
}{\:\!}$ we have 

$$ \frac{\dstl{\int_{\R^n}}|\nabla u(x)|^{p}\dfrac{1}{g(x)^{p\:\!-1}\;\!\!}dx}{\left(\dstl\int_{\R^n}|u(x)|^q v(x) g(x)\,dx\right)^{p/q}}\ge 
\frac{\dstl{\int_{\R^n}}|\nabla[{\;\!}{\cal R}_{\:\!\!g}^{}[u]](x)|^{p}\dfrac{1}{g(x)^{p\:\!-1}\;\!\!}dx}{\left(\dstl\int_{\R^n}{\cal R}_{\:\!\!g}^{}[u](x)^{q}v(x) g(x)\,dx\right)^{p/q}}.
$$

\end{prop}
\par\noindent{\bf Proof: }
Since $v$ is  non-increasing, we  have  from Proposition \ref{prop4.3} (iii) 
$
 {\cal R}_{\:\!\!g}^{}[v](x)=v(x)\;
 \mbox{ for  a.e. }x\in\R^{\;\!\!n}\bs{\;\!\!}\{0\}.
$
Hence, by Proposition \ref{prop4.4} (ii), we have 
$$\dstl{
 {\int_{\R^n}}|u(x)|^qv(x){\:\!}g(x){\:\!}dx
 \le {\int_{\R^n}}{\cal R}_{\:\!\!g}^{}[|u|^q](x){\;\!}{\cal R}_{\:\!\!g}^{}[v](x){\:\!}g(x){\:\!}dx=
  {\int_{\R^n}}{\cal R}_{\:\!\!g}^{}[u](x)^q{\;\!}v(x){\:\!}g(x){\:\!}dx.
}$$
Therefore the assertion holds  together with Proposition \ref{prop4.5}.\qed

\section{Preliminary results}
We begin to improve the  critical CKN-type inequality (\ref{CKNc}) {by using   poly-logarithm functions defined in  (\ref{poly-logarithm}).} Specifically we adopt   the poly-logarithmic  weight functions of  the form
\begin{equation}
|x|^{p-n}\left( \left( \log^k(R\eta/|x|)\right)^\a \; \prod_{j=1}^{k-1}\log^j(R\eta/|x|) \right)^{p-1} \qquad (\a\in \R, \; k\in \N)
\end{equation}
instead of $|x|^{p-n}$ on the left side of  (\ref{CKNc}). Here   we adopt the  notation  $\prod_{j=1}^{0}\log^j(R\eta/|x|)=1$.

Then we have  the following:

\begin{thm}\label{prop3.2} Let $n\in \N$, $1<p\le q <\infty$,  $\eta>0$ and $ 0 \le 1/p -1/q \le1/n$. 
Then,  we have  the followings:
\begin{enumerate}
\item
 If $\a\neq1$, $k\in \N$ and $R>\exp^{k-1}(1)$, then  there exists a positive  number $C_n$ such  that we  have 
 for any $u\in C_c^\infty(B_\eta\setminus \{0\})$ 
\begin{equation}
\begin{split}
\int_{B_\eta} &|\nabla u(x)|^p  |x|^{p-n}\left( \left( \log^k(R\eta/|x|)\right)^\a \; \prod_{j=1}^{k-1}\log^j(R\eta/|x|) \right)^{p-1}
\,dx\\ 
&\ge  C_n |\a-1|^{p-1+p/q} 
\left(  \int_{B_\eta} 
\frac{ |u(x)|^q \,dx}{|x|^n \dstl 
\left(\log^{k} (R\eta/|x|)\right)^{1+(1-\a)q/{p'}}\; \prod_{j=1}^{k-1}\log^j(R\eta/|x|)} \right)^{p/q}.
\end{split}\label{3.7}
\end{equation}
\item
If $\a=1$, $k\in \N \cup \{0\}$ and $R>\exp^{k}(1)$, then there exists a positive  number $C_n$ such  that  we  have 
 for any $u\in C_c^\infty(B_\eta \setminus \{0\})$ 
\begin{equation}
\begin{split}
\int_{B_\eta} &|\nabla u(x)|^p  |x|^{p-n}\left(  \prod_{j=1}^{k}\log^j(R\eta/|x|) \right)^{p-1}\,dx\\ 
&\ge  C_n \left(\dstl \int_{B_\eta} 
\frac{ |u(x)|^q \,dx}{|x|^n
 \dstl{\left(\log^{k+1} (R\eta/|x|)\right)^{1+q/{p'}}\prod_{j=1}^{k}\log^j(R\eta/|x|)}} \right)^{p/q}.
\end{split}\label{3.8}
\end{equation}
\item
Assume  that $C_n$ is the best constant. Then, $C_n$ satisfies the same assertions in Theorem \ref{th4.2} (ii).
\item 
If  $1<n< p$, $\a\le 1$  and $R$ is sufficiently large, then   
$C_n= C_{n,{\rm rad}}=S_{\rm rad}^{p,q;1/p'}(\R^n)$.
\end{enumerate}
\end{thm}
\medskip
\begin{rem}\label{remark3.1} \begin{enumerate}
\item 
The  inequality (\ref{3.8}) coincides with  (\ref{CKNc}) in the case $k=0$.
\item If $\a\le 1$, then
the weight functions in the inequalities (\ref{3.7}) and (\ref{3.8}) are locally integrable in $B_\eta$. 
Hence in this  case one can adopt $C^\infty_c(B_\eta) $ as a space of  test functions in place of $C^\infty_c(B_\eta\setminus \{0\}) $.
\item
If   $1<p \le (n+1)/2$, then   it  follows  from Theorem \ref{prop3.2} (iii),  Lemma \ref{lem2.1} and Remark \ref{remark2.2} (ii) that 
we  have (\ref{2.13}) for a  sufficiently large $R$, that is,
$C_n=C_{n,{\rm rad}}=  S^{p,q;1/p'}(\R^n)=S_{\rm rad}^{p,q;1/p'}(\R^n).$
\end{enumerate}
\end{rem}

\par\medskip

Theorem \ref{prop3.2}
  follows  from Theorem \ref{th4.2} with $w(t)=w_{k,\alpha}(t) \;(=w_{k,\alpha}(t;\eta))$ below.

\begin{defn}\label{feta}
Let $\a\in \R$, $\eta>0$,  $R>\exp^{k-1}(1)$ and $k\in \N$. 
\begin{enumerate}
\item
We  define $w_{k,\a}(t;\eta)$ such that $w_{k,\a}(t;\eta)\in C^{0,1}(\R_+)$, $w_{k,\a}(t;\eta)>0$ and 
\begin{equation}w_{k,\a}(t;\eta)=\begin{cases} t \left(\dstl\prod_{j=1}^{k-1}\log^{j}( R\eta/t)\right) \left(\log^{k}( R\eta/t)\right)^{\a}, &  0< t\le \eta,\\[4ex]
 \eta \left(\dstl\prod_{j=1}^{k-1}\log^{j} R\right) \left(\log^{k} R\right)^{\a}, & t\ge \eta.
\end{cases}
\end{equation}
Here  if $k=1$ we adopt the  notation  $\dstl\prod_{j=1}^{0} \log^j(R\eta/t)=1, \;\dstl\prod_{j=1}^{0}  \log^j R=1$.
\item
$w_{k,\a}(t;\eta)$ is abbreviated as $w_{k,\a}(t)$.
\end{enumerate}
\end{defn}
\par\medskip
Then we  see that  $w_{k,\a}(t)\in P(\R_+)$ if $ \a\le 1$, and $w_{k,\a}(t)\in Q(\R_+)$ if $ \a> 1$.\par\medskip
For $w(t)=w_{k,\a}(t)$  we have the followings:
\begin{enumerate}
\item When  $\a<1$, $R>\exp^k(1)$ and $\mu=\left(\log ^k R\right)^{1-\a}/(1-\a)$,
\begin{equation}
f_\eta (t)= \mu+\int_t^\eta\frac1{ w(s)}\,ds= \frac{1}{1-\a}
\left(\log ^k(R\eta/t)\right)^{1-\a}.
\end{equation}
\item When $\a=1$, $\mu=\log ^{k+1} R$ and $R>\exp^{k+1}(1)$,
\begin{equation}
f_\eta (t)=\mu+\int_t^\eta\frac1{ w(s)}\,ds= 
\log ^{k+1}(R\eta/t).\end{equation}
\item When  $\a>1$ and $R>\exp^k(1)$, 
\begin{equation}
f_\eta (t)=\int_0^t\frac1{ w(s)}\,ds= \frac{1}{\a-1}
\left(\log ^k(R\eta/t)\right)^{1-\a}.
\end{equation}

\end{enumerate}
Then  we   have 
{
\begin{align}& \frac{w(t)}{t}f_\eta(t)=\begin{cases}
 \dstl\frac{1}{1-\a}\prod_{j=1}^{k}\log ^j(R\eta/t),&\a<1,\\
 \dstl\prod_{j=1}^{k+1}\log ^j(R\eta/t),& \a=1,\\
 \dstl\frac{1}{\a-1}\prod_{j=1}^{k}\log ^j(R\eta/t),&\a>1,
\end{cases} \quad 0<t<\eta,.
\end{align}
}

Moreover  we  see that 
{
\begin{align}& \inf_{0<t\le \eta}  \frac{w(t)}{t}f_\eta(t)\ge \begin{cases}
 \dstl\frac{1}{1-\a}\prod_{j=1}^{k}\log ^j R,&\a<1,\\
 \dstl\prod_{j=1}^{k+1}\log ^j R,&\a=1,\\
 \dstl\frac{1}{\a-1}\prod_{j=1}^{k}\log ^j R,&\a>1.
\end{cases}
\end{align}
}
In particular if $R>\exp^{k}(1)$, then we have $C_0:= \dstl\inf_{0<t\le \eta}  \frac{w(t)}{t}f_\eta(t)>0$ so  that  $w(t)=w_{k,\a}(t)$ satisfies the non-degenerate condition (\ref{H-condition}). Since $ C_0\to \infty$ as $R\to\infty$, $C_0\ge 1 $ holds  for a sufficiently large $R>\exp^{k}(1)$.
\medskip
\par\noindent{\bf {Proof of  Theorem \ref{prop3.2}:} } The  assertions  (i), (ii) and  (iii) follow from Theorem \ref{th4.2}. Hence, 
it suffices to  show the assertion (iv). 
From Remark \ref{remark3.1} (ii),  we adopt $C^\infty_c(B_\eta)$ as the space of  test functions.
We use Proposition \ref{prop4.6} with appropriate $g$ and $v$.  To  this  end 
we define $g(x)=g(|x|)$ by
\begin{equation}g(t)^{1-p} = w_{k,\a}(t)^{p-1}t^{1-n}\qquad (t= |x|, \; t\le \eta),\label{defofg}\end{equation}
that is
$$ g(t)= \frac{t^{(n-p)/(p-1)} }{\left(\dstl\prod_{j=1}^{k-1}\log^{j}( R\eta/t)\right)\left(\log^{k}(R\eta/t)\right)^{\a}}\qquad ( 1<n<p,\; {\a\le 1}).$$
We need to  check that  $g(t)$ is decreasing for a sufficiently large $R$. By a direct calculation we  have 
\begin{equation}\begin{split}
g'(t)&=  \frac{{t^{(n-p)/(p-1)-1}}}{\left(\dstl\prod_{j=1}^{k-1}\log^{j}( R\eta/t)\right)\left(\log^{k}(R\eta/t)\right)^{\a}}
\left(\frac{n-p}{p-1}+ \sum_{j=1}^{k-1} \frac{1}{\dstl\prod_{l=1}^j \log^l(R\eta/t)}
+ \frac{\a}{\dstl\prod_{l=1}^k \log^l(R\eta/t)}\right)\\
& \le  \frac{{t^{(n-p)/(p-1)-1}}}{\left(\dstl\prod_{j=1}^{k-1}\log^{j}( R\eta/t)\right) \left(\log^{k}(R\eta/t)\right)^{\a}}
\left(\frac{n-p}{p-1}+ \sum_{j=1}^{k-1} \frac{1}{\dstl\prod_{l=1}^j \log^l R}
+ \frac{\a}{\dstl\prod_{l=1}^k \log^l R}\right).
\end{split}
\end{equation} Since $n<p$, we have 
  $g'(t)\le 0$ provided that $R$ is sufficiently  large.
We define   $v(x)=v(|x|)$ for $t=|x|$ by

\begin{equation}
 v(t)=\begin{cases} t ^{-p'(n-1)}\left(\log ^k(R\eta/t)\right)^{-(1-\a)(1+q/p')} & (\a<1),\\[2ex]
 t ^{-p'(n-1)}\left(\log ^{k+1}(R\eta/t)\right)^{-(1+q/p')}& (\a=1).
 \end{cases}\label{defofv}
\end{equation}
Then we  see that 
\begin{equation} v(t)g(t)=\begin{cases}\dstl
{ t^{-n}\left(\dstl{\prod_{j=1}^{k-1}\log^j (R\eta/t)}\right)^{-1}
 \dstl{\left(\log^{k} (R\eta/t)\right)^{-1-(1-\a)q/{p'}}}}  &(\a<1),\\
\dstl { t^{-n}\left(\dstl{\prod_{j=1}^{k}\log^j (R\eta/t)}\right)^{-1}
 \dstl{\left(\log^{k+1} (R\eta/t)\right)^{-1-q/{p'}}}} &(\a=1).
\end{cases}
\end{equation}
In order to apply Proposition \ref{prop4.6},
we need to  check if $ v'(t)\le 0$.  
By a direct calculation we  have

\begin{equation}
v'(t)= \begin{cases} 
-p'(n-1) t^{-p'(n-1)-1}(\log^k(R\eta/t))^{- (1-\a)(1+ q/p')-1} &\\
\quad \times\left( \log^k(R\eta/t) -\frac{(1-\a)(1+q/p')}{ p'(n-1) }\frac{1}{\prod_{j=1}^{k-1}\log^j(R\eta/t)}\right)&{ (\a<1)},\\[4ex]
 -p'(n-1) t^{-p'(n-1)-1}(\log^{k+1}(R\eta/t))^{-1} &\\
\quad\times \left( \log^{k+1}(R\eta/t) -\frac{1+q/p'}{ p'(n-1) }\frac{1}{\prod_{j=1}^{k}\log^j(R\eta/t)}\right)&  (\a=1).
 \end{cases}
\end{equation}

First  we assume that $\a\neq 1$. Then  we  show 
$$ \log^k(R\eta/t) \ge \frac{(1-\a)(1+q/p')}{ p'(n-1)} \frac{1}{\dstl \prod_{j=1}^{k-1}\log^j(R\eta/t)} \qquad (0<t\le \eta). $$
Since $\log^j(R\eta/t)\; ( j=1,2,\cdots, k)$ is decreasing, it  suffices  to  check  
  $$\log^k R \ge \frac{(1-\a)(1+q/p')}{ p'(n-1)} \frac{1}{\dstl\prod_{j=1}^{k-1}\log^j R},$$ and  this    is  valid for a   sufficiently large $R$.
Secondly we  assume that $\a=1$. Then, in a similar way we see that
  $$ \log^{k+1} R \ge \frac{1+q/p'}{ p'(n-1)} \frac{1}{\dstl\prod_{j=1}^{k}\log^j R}$$
  holds for a sufficiently large $R$, hence  the assertion follows.
 \par\noindent{\bf End of  the proof:} By Proposition \ref{prop4.6}  with  (\ref{defofg}) and (\ref{defofv}), 
 we can  assume  that 
$u$  is radial. 
Then we see that test function $u$ for Theorem \ref{prop3.2} can be assumed to be radial.
The assertion (iv) therefore  comes from Remark \ref{remark2.2} (ii). \qed

\section{Main results}
We improve  Theorem \ref{prop3.2} by using the   super-logarithm
 defined in  \S 2.1 as  follows.
Specifically we adopt   the super-logarithmic  weight functions of  the form
\begin{equation}
|x|^{p-n} \left(  \left( \prod_{l=1}^\infty \frac{ F^l\left(\frac{a\eta}{|x|}\right)}{a}
\right)
\left( F^k \left(L\left(\frac{\eta}{|x|}\right)+a\right) \right)^\a \;
\prod_{j=0}^{k-1} F^j\left(L\left(\frac{\eta}{|x|}\right)+a\right)
\right)^{p-1}\label{SL}
\end{equation}
instead of 
$$
|x|^{p-n}\left( \left( \log^k(R\eta/|x|)\right)^\a \; \prod_{j=1}^{k-1}\log^j(R\eta/|x|) \right)^{p-1} \qquad (\a\in \R, \; k\in \N)
$$
on the left side of  (\ref{3.7}). 
Then we have  the following:

\begin{thm}\label{thm7.1}
Let $a>1, k\in \N \cup \{0\}$, $\a\in \R$ and $\eta>0$. 
Assume  that  $n\in \N$, $1<p\le q <\infty$ and  $ 0 \le 1/p -1/q \le1/n$ . 
Then, we have  the followings:
\par\noindent (i)  When $\a\neq1$, there exists a positive number $C_n$ such  that
 we  have 
 for any $u\in C_c^\infty(B_\eta\setminus \{0\})$ 
\begin{equation}
\begin{split}
&\int_{B_\eta} |\nabla u(x)|^p \;
|x|^{p-n} \left(  \left( \prod_{l=1}^\infty \frac{ F^l\left(\frac{a\eta}{|x|}\right)}{a}
\right)
\left( F^k \left(L\left(\frac{\eta}{|x|}\right)+a\right) \right)^\a \;
\prod_{j=0}^{k-1} F^j\left(L\left(\frac{\eta}{|x|}\right)+a\right)
\right)^{p-1}\,dx\\ 
&\\
&\ge  C_n |\a-1|^{p-1+p/q}\times \\
&\left( \int_{B_\eta} 
\frac{ |u(x)|^q \,dx}{ |x|^n\left( \dstl\prod_{l=1}^\infty \frac{ F^l\left(\frac{a\eta}{|x|}\right)}{a}
\right)\left( F^k\left(L\left(\frac{\eta}{|x|}\right)+a\right)\right)^{1+(1-\a)q/p'} 
\dstl \prod_{j=0}^{k-1}F^j\left(L\left(\frac{\eta}{|x|}\right)+a\right)}\right) ^{p/q},
\end{split}\label{4.11}
\end{equation}
where  if $k=0$ we adopt the  notation  $\dstl\prod_{j=0}^{k-1} F^j(u)=1$.\par
Here 
$
 F(u)=F_a(u)=a-\log a+\log u \; (u\ge a), \;  F^k(u)= F(F^{k-1}(u))$ for $ u\ge a >1$  and $L$ is the super-logarithm defined by 
 Definition \ref{defn:slog}.

\par\noindent (ii)  When $\a=1$, there exists a positive number $C_n$ such  that we  have 
 for any $u\in C_c^\infty(B_\eta\setminus \{0\})$ 
 
\begin{equation}
\begin{split}
&\int_{B_\eta} |\nabla u(x)|^p |x|^{p-n} \left(  \left( \prod_{l=1}^\infty \frac{ F^l\left(\frac{a\eta}{|x|}\right)}{a}
\right)
\prod_{j=0}^{k} F^j\left(L\left(\frac{\eta}{|x|}\right)+a\right)
\right)^{p-1}\,dx\\ 
&\ge  C_n \left( \int_{B_\eta} 
\frac{ |u(x)|^q \,dx}{ |x|^n\left( \dstl\prod_{l=1}^\infty \frac{ F^l\left(\frac{a\eta}{|x|}\right)}{a}
\right)\left( F^{k+1}\left(L\left(\frac{\eta}{|x|}\right)+a\right)\right)^{1+q/p'} 
\dstl \prod_{j=0}^{k}F^j\left(L\left(\frac{\eta}{|x|}\right)+a\right)} \right)^{p/q}.
\end{split}\label{4.12}
\end{equation}

\par\noindent (iii)
Assume  that $C_n$ is the best constant. Then, $C_n$ satisfies the same assertions in Theorem \ref{th4.2} (ii).
Further if  $ a> 1$  and $a\ge  |\a-1|^{1/(k+1)}$, then we  have  $C_0\ge 1$.
\par\noindent (iv)
If  $1<n< p$, $\a\le 1$  and $a$ is sufficiently large, then   
$C_n= C_{n,{\rm rad}}=S_{\rm rad}^{p,q;1/p'}(\R^n)$.

\end{thm}

\begin{rem}\label{remark4.1}
 \begin{enumerate}\item 
The  inequalities (\ref{4.11}) and (\ref{4.12}) improve (\ref{CKNc}), (\ref{3.7}) and (\ref{3.8}).
In particular if $ k=0$ in (\ref{4.12}), then the weight in  the left-hand side is simply
$$ |x|^{p-n}  \left( \prod_{l=1}^\infty \frac{ F^l\left({a\eta}/{|x|}\right)}{a}
\right)^{p-1}.$$
Roughly speaking this can be  thought of as a formal  limit of 
$$|x|^{p-n}\left(  \prod_{l=1}^{k}\log^l(R\eta/|x|) \right)^{p-1}$$
as $k\to\infty$.
\item If $\a\le 1$, then
the weight functions in the inequalities (\ref{4.11}) and (\ref{4.12}) are locally integrable in $B_\eta$. Hence one can adopt $C^\infty_c(B_\eta) $ as a space of  test functions in place of $C^\infty_c(B_\eta\setminus \{0\}) $.
\item
If   $1<p \le (n+1)/2$, then   it  follows  from Theorem \ref{thm7.1} (iii),  Lemma \ref{lem2.1} and Remark \ref{remark2.2} (ii) that 
we  have (\ref{2.13}) for $ a> 1$  and $a\ge  |\a-1|^{1/(k+1)}$, that is,
$C_n= C_{n,{\rm rad}}= S^{p,q;1/p'}(\R^n)=S_{\rm rad}^{p,q;1/p'}(\R^n).$
\end{enumerate}
\end{rem}

Theorem \ref{thm7.1}will be shown  by using  Theorem \ref{th4.2} with $w(t)=w^{sl}_{k,\a}(t)\; (=w^{sl}_{k,\a}(t;\eta))$
below.
\begin{defn}
\label{Definition4.1}
\begin{enumerate}
\item
We define $w^{sl}_{k,\a}(t;\eta)\in C^{0,1}(\R_+)$ for $k\in  \N \cup \{0\}$ such  that
\begin{equation}\begin{split}
w^{sl}_{k,\a}&(t;\eta)=\\
&\begin{cases}t\left( \dstl\prod_{l=1}^\infty \frac{ F^l\left({a\eta}/{t}\right)}{a}
\right) \left(\dstl\prod_{j=0}^{k-1} F^j\left(L\left({\eta}/{t}\right)+a\right)\right)\; \left(F^k\left(L\left({\eta}/{t}\right)+a\right)\right)^\a 
&(0<t\le \eta),\\[4ex]
 \eta a^{k+\a} & (t\ge \eta). \end{cases}\label{4.14}
\end{split}
\end{equation}
Here  if $k=0$ we adopt the  notation  $\dstl\prod_{j=0}^{-1} F^j\left(L\left({\eta}/{t}\right)+a\right)=1$.
\item  $w^{sl}_{k,\a}(t;\eta)$ is  abbreviated as $w^{sl}_{k,\a}(t)$.
\end{enumerate}
\end{defn}
Then we  have  the following:

\begin{lem} Let $w(t)=w^{sl}_{k,\a}(t)\in W(\R_+)$ be given by (\ref{4.14}). Then,
\begin{enumerate}
\item   The  corresponding $f_\eta(t)\;  (0<t<\eta)$  is given by
\begin{equation}f_\eta(t)=\begin{cases} \mu+ \int_t^\eta \frac{1}{w(s)}\,ds= 
 \frac{1}{1-\a} \left(F^k\left(L\left({\eta}/{t}\right)+a\right)\right)^{1-\a}
& (\a<1, \mu=\frac{1}{1-\a}a^{1-\a}), \\ 
 \mu+\int_t^\eta \frac{1}{w(s)}\,ds= F^{k+1}\left(L\left({\eta}/{t}\right)+a\right)&(\a=1, \mu= a),\\
\int_0^t\frac{1}{w(s)}\,ds= \frac{1}{\a-1}\left( F^k\left(L\left({\eta}/{t}\right)+a\right)\right)^{1-\a} &(\a>1).
\end{cases} \label{4.19}
\end{equation}
\item $\dstl\lim_{t\to+0}\frac{w(t)}{t}f_\eta(t)=\infty$. In particular  if $ a> 1$ and $a\ge |\a-1|^{1/(k+1)}$, then  $w(t)$ satisfies the non-degenerate condition (\ref{H-condition}) with $C_0\ge 1$.
\end{enumerate}
\end{lem} 
\par\noindent{\bf Proof of (i): }  Let $w(t)=w^{sl}_{k,\a}(t)$. First  we assume  that $\a< 1$.
Since it  holds  that {by (\ref{2.34})}
\begin{equation} \frac{1}{1-\a}\frac{d}{dt}\left(F^k\left(L\left(t\right)+a\right)\right)^{1-\a}= 
\frac{1}{t\left( \dstl\prod_{l=1}^\infty \frac{ F^l\left(at\right)}{a}
\right)\left( \dstl\prod_{j=0}^{k-1} F^j\left(L\left(t\right)+a\right)\right)\;\left(F^k\left(L\left(t\right)+a\right)\right)^\a},\label{4.15}\end{equation}
 we see that
\begin{equation}
\begin{split} \int_t^\eta\frac{1}{w(s)}&\,ds
=\int^{\eta/t}_1 \frac{1}{u\,\left( \dstl\prod_{l=1}^\infty \frac{ F^l\left(au\right)}{a}
\right)\left( \dstl\prod_{j=0}^{k-1} F^j\left(L\left(u\right)+a\right)\right)\; \left(F^k\left(L\left(u\right)+a\right)\right)^\a}\,du 
 \quad(\eta/s=u)\\
&=\int^{\eta/t}_1 \frac{1}{1-\a}\frac{d}{du}\left( F^k\left(L\left(u\right)+a\right)\right)^{1-\a}\,du\\
&= \frac{1}{1-\a}\left( \left( F^k\left(L\left({\eta}/{t}\right)+a\right)\right)^{1-\a}-\left(F^k\left(L\left(1\right)+a\right)\right)^{1-\a}\right)\\
&=\frac{1}{1-\a}\left( F^k\left(L\left({\eta}/{t}\right)+a\right)\right)^{1-\a}-\mu
\quad\left(\mu=\frac{1}{1-\a}\left( F^k\left(L\left(1\right)+a\right)\right)^{1-\a} = \frac{a^{1-\a}}{1-\a}\right).
\end{split}\label{4.16}
\end{equation}
Secondly  we  assume  that $\a>1$.
By Theorem \ref{thm:F} (ii) we  have  $\dstl \lim_{t\to+0}F^k\left(L\left({\eta}/{t}\right)+a\right)=\infty$. 
Then in a similar way  we have
\begin{equation}
\begin{split} \int_0^t\frac{1}{w(s)}\,ds
&= \frac{1}{\a-1}\left( F^k\left(L\left({\eta}/{t}\right)+a\right)\right)^{1-\a}.
\end{split}\label{4.16'}
\end{equation}
Thirdly we assume that $\a=1$. 
Noting  that $$F^{k+1}\left(L\left(t\right)+a\right)=F(F^k\left(L\left(t\right)+a\right))
\;\text{ and }\;\dstl \frac{d}{dt}F^{k+1}\left(L\left({\eta}/{t}\right)+a\right)
= \frac{1}{t\left( \dstl\prod_{l=1}^\infty \frac{ F^l\left(at\right)}{a}
\right) \dstl\prod_{j=0}^{k} F^j\left(L\left(t\right)+a\right)},$$
we  have
\begin{equation}
\begin{split} \int_t^\eta\frac{1}{w(s)}\,ds
&  =F^{k+1}\left(L\left({\eta}/{t}\right)+a\right)-\mu \qquad ( w(s)=w^{sl}_{k,1}(s),\; \mu=F^{k+1}\left(L\left(1\right)+a\right) =a).
\end{split}\label{4.18}
\end{equation}
Therefore
we  see that if $\a\le1$, $w^{sl}_{k,\a}(t)\in P(\R_+)$ and  if $\a>1$,  $w^{sl}_{k,\a}(t)\in  Q(\R_+)$.
By Definition \ref{Definition2.3}   we  immediately have (\ref{4.19}) from (\ref{4.16}), (\ref{4.16'}) and (\ref{4.18}). \qed
\par\medskip\noindent
{\bf Proof of (ii).} Let $w(t)=w^{sl}_{k,\a}(t)\in W(\R_+)$.
From 
 (\ref{4.19}) we  have
\begin{equation}
\frac{w(t)}{t}f_\eta(t)
=\begin{cases} 
\frac{1}{1-\a}\left( \dstl\prod_{l=1}^\infty \frac{ F^l\left(a\eta/t\right)}{a}
\right)\dstl \prod_{j=0}^kF^j\left(L\left({\eta}/{t}\right)+a\right)&(\a<1),\\[4ex]
\left( \dstl\prod_{l=1}^\infty \frac{ F^l\left(a\eta/t\right)}{a}
\right)\; \left(\dstl \prod_{j=0}^kF^j\left(L\left({\eta}/{t}\right)+a\right)\right)F^{k+1}\left(L\left({\eta}/{t}\right)+a\right)&  (\a=1),\\[4ex]
\frac{1}{\a-1}\left( \dstl\prod_{l=1}^\infty \frac{ F^l\left(a\eta/t\right)}{a}
\right)\dstl \prod_{j=0}^kF^j\left(L\left({\eta}/{t}\right)+a\right)&(\a>1).

\end{cases}
\end{equation}
Because  $\prod_{l=1}^\infty \frac{ F^l\left(a\eta/t\right)}{a}
\ge 1$, $\dstl \lim_{t\to+0}F^k(L(\eta/t)+a)=\infty$ and $ F^k\left(L\left(\eta/t\right)+a\right)\ge a$  for $k\in \N\cup \{0\}$ and  $0<t<\eta$, 
we  have $\dstl\lim_{t\to+0} \frac{w(t)}{t}f_\eta(t)=\infty$ and the estimate
 \begin{equation}
\frac{w(t)}{t}f_\eta(t)
\ge\begin{cases} 
\frac{a^{k+1}}{|1-\a|}& (\a\neq1),\\
 a^{k+2}& (\a=1).
\end{cases}\label{4.10}
\end{equation}
Therefore,   the assertion (ii) is satisfied, if $ a>1$  and $a\ge  |\a-1|^{1/(k+1)}$.
\qed

\medskip
\par\noindent{\bf Proof of Theorem  \ref{thm7.1}: }
It suffices to  show the assertion (iv). 
From Remark \ref{remark4.1} (ii),  we adopt $C^\infty_c(B_\eta)$ as the space of  test functions.
As before we apply Proposition \ref{prop4.6} with appropriate $g$ and $v$. 
Then we see that test function $u\in C^\infty_c(B_\eta)$ for Theorem \ref{thm7.1} can be assumed to be radial. The assertion (iv) therefore  comes from Remark \ref{remark2.2} (ii).
\par
To  this  end 
we define $g(x)=g(|x|)$ by
\begin{equation}\label{g(t)} g(t)^{1-p} = w^{sl}_{k,\a}(t)^{p-1}t^{1-n}\qquad (t= |x|, \; t\le \eta).
\end{equation}
Then 
$$ g(t)= \frac{t^{({n-p})/({p-1})} }{\left( \dstl\prod_{l=1}^\infty \frac{ F^l\left(a\eta/t\right)}{a}
\right) \left(\dstl\prod_{j=0}^{k-1} F^j\left(L\left({\eta}/{t}\right)+a\right)\right)\; \left(F^k\left(L\left({\eta}/{t}\right)+a\right)\right)^\a}\qquad\qquad ( 
\a\le 1).
$$

We  also define  $ v_k(x)=v_k(|x|)$ for  $k\in \N\cup\{0\}$ and $t=|x|$ by 
 \begin{equation}
 v_k(t)=\begin{cases} t ^{-p'(n-1)}\left(F^k\left(L\left({\eta}/{t}\right)+a\right)\right)^{-(1-\a)(1+q/p')} & (\a<1),\\
 t ^{-p'(n-1)}\left(F^{k+1}\left(L\left({\eta}/{t}\right)+a\right)\right)^{-(1+q/p')}&(\a=1).
 \end{cases}\label{decreasing}
\end{equation}
Then we show  the following:

\begin{lem}\label{lemma4.3} Let $k\in \N\cup \{0\}$.
\begin{enumerate}
 \item
 If  $a$ is sufficiently large, then $g(t)\; (0<t\le \eta) $ is decreasing.
\item
 If $a$  is sufficiently large, then $v_k(t)\; (0<t\le \eta) $ is decreasing. 
 
 \end{enumerate}
\end{lem}
\par\noindent{\bf Proof of (i):}  We  note  the following relations:
For any $ k \in \N $ and  $a>1$, 
\begin{equation}
\begin{cases}
\dstl  L'(r)=\frac{d}{{dr}}(a+L(r))= \frac{a}{ \tF(ar)}=\frac{1}{{r}\displaystyle\prod_{j=1}^{\infty}\frac{F^{j}(ar)}{a}}\qquad \qquad &(r>1),\\
\dstl \frac{d}{dr}F^k(a+L(r))
= \dstl\frac{L'(r)}{ \dstl\prod_{j=0}^{k-1} F^j(a+L(r))} &(r>1),\\[2ex]
\dstl \frac d{du}\tF(u)= \tF(u) \left(\dstl\sum_{k=0}^{\infty}\frac1{\prod_{j=0}^{k}F^{j}(u)}\right)&(u>a),
\\[4ex]
\dstl
\frac{d}{dr}\left(\displaystyle\prod_{j=1}^{\infty}\frac{F^{j}(ar)}{a}\right)=
\left(\displaystyle\prod_{j=1}^{\infty}\frac{F^{j}(ar)}{a}\right)\sum_{k=1}^\infty\frac{1}{{r} \displaystyle\prod_{j=1}^k F^j(ar)} & (r>1).
\end{cases}\label{4.15}
\end{equation}

By (\ref{4.15}) we have \begin{equation}
g'(t)= \frac{g(t)}{t}  G(t),\qquad{ (0<t<\eta)}\end{equation}
where
 \begin{equation}\begin{split}
 G(t)&= \frac{n-p}{p-1}+ \sum_{k=1}^\infty \frac{1}{ \dstl\prod_{j=1}^k F^j(a\eta/t)}
+\sum_{j=0}^{k-1} \frac{1}{ \left( \dstl\prod_{l=1}^\infty \frac{ F^l\left(a\eta/t\right)}{a}
\right)\dstl\prod_{i=0}^j F^i\left(L\left({\eta}/{t}\right)+a\right)}\\
&
+ \a \frac{1}{ \left( \dstl\prod_{l=1}^\infty \frac{ F^l\left(a\eta/t\right)}{a}
\right)\dstl\prod_{i=0}^k F^i\left(L\left({\eta}/{t}\right)+a\right)}.
\end{split}
\end{equation}
Then $G(t)$ satisfies 
$$ G(t)\le  \frac{n-p}{p-1}+ 
\sum_{k=1}^\infty \frac{1}{ a^k}
+\sum_{j=0}^{k-1} \frac{1}{a^{j+1}}
+|\a |\frac{1}{a^{k+1}} <\frac{n-p}{p-1} + \frac{2}{a-1} + |\a|\frac{1}{a^{k+1}}.
$$
Since $p>n$, we  see that  $g'(t)<0 $ if $a$ is  sufficiently large. \qed
\par\medskip
\noindent{\bf Proof of (ii): }
 First  we assume  that $\a=1$. Then
\begin{equation} \begin{split}
v'_k(t)&= -p'(n-1) t^{-p'(n-1)-1} \left( F^{k+1}\left(L\left({\eta}/{t}\right)+a\right)\right)^{-1+ q/p'-1}\\
&\times\left( F^{k+1}\left(L\left({\eta}/{t}\right)+a\right) - \frac {1+ q/p'}{p'(n-1)}\frac{1}{\left( \dstl\prod_{l=1}^\infty \frac{ F^l\left(a\eta/t\right)}{a}
\right)\dstl\prod_{j=0}^kF^j\left(L\left({\eta}/{t}\right)+a\right)}\right).
\end{split}
\end{equation}
It suffices to  show that if  $a$ is  large,  then
\begin{equation} 
 F^{k+1}\left(L\left({\eta}/{t}\right)+a\right) \ge  \frac {1+ q/p'}{p'(n-1)}\frac{1}{\left( \dstl\prod_{l=1}^\infty \frac{ F^l\left(a\eta/t\right)}{a}
\right)\dstl\prod_{j=0}^kF^j\left(L\left({\eta}/{t}\right)+a\right)}.
\end{equation}
Since $F^j\left(L\left({\eta}/{t}\right)+a\right)\; ( j\in \N\cup\{0\}) $ and $ \prod_{l=1}^\infty \frac{ F^l\left(a\eta/t\right)}{a}
$ are decreasing,  we  have 
$$  F^{k+1}\left(L\left({\eta}/{t}\right)+a\right)\ge F^{k+1}\left(L\left(1\right)+a\right) = a $$
and
$$\frac{1}{\left( \dstl\prod_{l=1}^\infty \frac{ F^l\left(a\eta/t\right)}{a}
\right)\dstl\prod_{j=0}^kF^j\left(L\left({\eta}/{t}\right)+a\right)} \le \frac{1}{\prod_{j=0}^kF^j\left(L\left(1\right)+a\right)}= \frac 1{a^{k+1}}.$$
Hence, if $$ a^{k+2} \ge(1+q/p')/ (p'(n-1))$$ is  satisfied, then $v_k'(t)< 0\; (0<t{< }\eta)$. Therefore we have the  desired result.
\par
 We proceed to the  case that $\a\neq 1$.  Then we have 
\begin{equation} \begin{split}
v_k'(t)  &= -p'(n-1) t^{-p'(n-1)-1}\left(F^k\left(L\left({\eta}/{t}\right)+a\right)\right)^{- (1-\a)(1+ q/p')-1}\\
&\times
\left( F^k\left(L\left({\eta}/{t}\right)+a\right) - \frac { (1-\a)(1+ q/p')}{ p'(n-1)}\frac{1}{\left( \dstl\prod_{l=1}^\infty \frac{ F^l\left(a\eta/t\right)}{a}
\right)\dstl\prod_{j=0}^{k-1}F^j\left(L\left({\eta}/{t}\right)+a\right)}\right),\end{split}
\end{equation}
where  we use  the  notation  $\prod_{j=0}^{-1}F^j\left(L\left({\eta}/{t}\right)+a\right)=1$. 
In a similar way,  if $$ a^{k+1}\ge  (1-\a) (1+q/p')/ (p'(n-1))$$ holds, then we  see that $v_k'(t)< 0\; {(0<t<\eta)}$, and the assertion is proved.\qed

\par\medskip\noindent{\bf End of the proof of Theorem  \ref{thm7.1} (iv): }
We  have
\begin{equation*} \begin{split}&
v_k(t)g(t)=\\
&\begin{cases}\dstl
 t^{-n} \left(\left( \dstl\prod_{l=1}^\infty \frac{ F^l\left(a\eta/t\right)}{a}
\right)  \left(\dstl\prod_{j=0}^{k-1}F^j\left(L\left({\eta}/{t}\right)+a\right)\right)\right)^{-1}
 \dstl{\left(F^k\left(L\left({\eta}/{t}\right)+a\right)\right)^{-1-(1-\a)q/{p'}}} &{(\a<1)},\\
 t^{-n}\left(\left( \dstl\prod_{l=1}^\infty \frac{ F^l\left(a\eta/t\right)}{a}
\right)
 \left(\dstl\prod_{j=0}^{k}F^j\left(L\left({\eta}/{t}\right)+a\right)\right)\right)^{-1}\left(F^{k+1}\left(L\left({\eta}/{t}\right)+a\right)\right)^{-(1+q/p')}  &(\a=1).
\end{cases}
\end{split}
\end{equation*}
Then,
the  best constant $C_n$ is given by
\begin{align*}
C_n&=  \inf_{u\in C^\infty_c(B_\eta)\setminus \{0\}}\frac{\dstl{\int_{B_\eta}}|\nabla u(x)|^{p}w_{k,\a}^{sl}(|x|)^{p-1}|x|^{1-n}
dx}{\left(\dstl\int_{B_\eta}\dfrac{|u(x)|^q|x|^{1-n} }{w_{k,\a}^{sl}(|x|) f_\eta(|x|)^{1+q/p'}}\,dx\right)^{p/q}}
=\inf_{u\in C^\infty_c(B_\eta)\setminus \{0\}}\frac{\dstl{\int_{B_\eta}}|\nabla u(x)|^{p}\dfrac{1}{g(x)^{p\:\!-1}\;\!\!}dx}{\left(\dstl\int_{B_\eta}|u(x)|^q v_k(x) g(x)\,dx\right)^{p/q}}.
\end{align*}
As   functions $ g(t)$ and  $v(t)$ in Proposition \ref{prop4.6}, we adopt $g(t)$ defined by (\ref{g(t)}) and $v_k(t) $ defined by (\ref{decreasing}) respectively. Noting that  $C^\infty_c(B_\eta)$ can be  replaced by $C_c^{0,1}(B_\eta)$ using  a density argument, 
it follows from 
 Proposition \ref{prop4.6} that
 $u$  can be assumed   radial. 
 Then we can adopt $C_c^\infty(B_\eta)_{\rm rad}$ as a space  of test functions
 in Theorem \ref{thm7.1}.
The assertion (iv) therefore  comes from Remark \ref{remark2.2} (ii). \qed
\section{Further results}
If  $p=q$ holds,  then the inequalities  (\ref{4.11}) and (\ref{4.12}) become the Hardy-type.
In this  case we can further improve  them by adding  sharp remainders. 
To  this end  we recall and modify the one-dimensional Hardy's inequalities in \cite{ho3}.
The following is an easy corollary to Theorem 3.1 in \cite{ho3}.
\begin{prop} \label{T3.1}\par \noindent
Assume that $1<p<\infty$, $\eta>0$, $\mu>0$ and  $w\in W(\mathbf R_+)$. Then,
there exist  a positive number  $C=C(w,p,\eta,\mu )$  such  that  for  every 
$u\in C^\infty_c((0,\eta))$, we  have 
\begin{equation}
\begin{split}  \int_0^\eta   |u'(t)|^p w(t)^{p-1}\,dt  
\ge &  (1/p')^p
 \int_0^\eta \frac{ |{u(t)}|^p  \,dt }{w(t) f_\eta(t)^p} + C\int_0^\eta \frac{ |u(t)|^p  \,dt}{ w(t) f_\eta(t)^p G_\eta(t)^2},\label{nct2}
\end{split}
\end{equation}
where $f_\eta(t)$ and $G_\eta(t)$ are  defined by Definition \ref{Definition2.3} and the coefficient $ (1/p')^p$ is  sharp.
\end{prop}
Then we immediately have the $n$-dimensional version:

\begin{cor} \label{cor5.2}\par \noindent
Assume that $n\in \N$,  $1<p<\infty$, $\eta>0$, $\mu>0$ and  $w\in W(\mathbf R_+)$. Then,
there exist  a positive number  $C=C(w,p,\eta,\mu )$  such  that  for  every 
$u\in C^\infty_c(B_\eta \setminus \{0\})$, we  have 
\begin{equation}
\begin{split}  \int_{B_\eta}   |\nabla u(x)|^p w(|x|)^{p-1}|x|^{1-n}\,dx  
\ge &  (1/p')^p
 \int_{B_\eta} \frac{ |{u(x)}|^p |x|^{1-n} \,dx }{w(|x|) f_\eta(|x|)^p} + C\int_{B_\eta} \frac{ |{u(x)}|^p |x|^{1-n} \,dx}{ w(|x|) f_\eta(|x|)^p G_\eta(|x|)^2}.\label{nct2}
\end{split}
\end{equation}
\end{cor}

From Corollary \ref{cor5.2} with $w(t)= w^{sl}_{k,1}(t)$ we  have the following.
For the sake of simplicity, the results are limited to case that $\a=1$.
By  Remark \ref{remark4.1} (ii) we can 
 adopt  $C^\infty_c(B_\eta)$ as  the space of test functions. (See also Remark \ref{remark3.1}(ii)). 
\begin{prop}\label{Proposition5.1}Let $k\in \N \cup \{0\}$ and $\eta>0$. 
Let $w^{sl}_{k,1} \in W(\R_+)$  be given by Definition \ref{Definition4.1}.
Assume  that  $n\in \N$, $1<p$ and $ a> 1$. 
{ Then, there exists a positive number $C=C(k,p,\eta,a)$  such  that
 for any $u\in C^\infty_c(B_\eta)$}

\begin{equation}
\begin{split}
&\int_{B_\eta} |\nabla u(x)|^p  w^{sl}_{k,1}(|x|)^{p-1}  |x|^{1-n}\,dx
\ge  (1/p')^p\int_{B_\eta} 
\frac{ |u(x)|^p |x|^{1-n}\,dx}{ w^{sl}_{k,1}(|x|)\left(F^{k+1}\left(L\left({\eta}/{|x|}\right)+a\right)\right)^{p}} \\ &\quad + 
C \int_{B_\eta} 
\frac{ |u(x)|^p |x|^{1-n}\,dx}{ w^{sl}_{k,1}(|x|)\left(F^{k+1}\left(L\left({\eta}/{|x|}\right)+a\right)\right)^{p}\left(a-\log a
+ \log \left(F^{k+1}\left(L\left({\eta}/{|x|}\right)+a\right) \right)\right)^2} .
\end{split}\label{5.6}
\end{equation}
Here,  $w^{sl}_{k,1}(t)$  is  defined by (\ref{4.14}) and the coefficient $ (1/p')^p$ is  sharp.
\end{prop}
\noindent{\bf Proof:} By virtue  of  (\ref{fetaP})  and (\ref{Geta}) in Definition \ref{Definition2.3}, 
we can calculate $G_\eta(t)$ to  obtain 
\begin{equation} G_\eta(t)=\mu -\log\mu + \log \left(\mu+ \int_t^\eta\frac1{w(s)}\,ds\right).
\end{equation}
 By  (\ref{4.18}) we have  
 \begin{equation} G_\eta(t; w^{sl}_{k,1},\mu)=\mu -\log\mu + \log \left(F^{k+1}\left(L\left({\eta}/{t}\right)+a\right)\right),
\end{equation}
where
$\mu=F^{k+1}\left(L\left(1\right)+a\right){=a}$. Therefore  we  have  (\ref{5.6}).\qed

\bigskip\noindent
{\large\bf Hiroshi Ando, Toshio Horiuchi, Eiichi Nakai \\Department of Mathematics\\Faculty of Science \\ Ibaraki University\\
Mito, Ibaraki, 310, Japan}\par\bigskip\noindent
e-mail: \par\noindent
hiroshi.ando.math@vc.ibaraki.ac.jp\par\noindent
toshio.horiuchi.math@vc.ibaraki.ac.jp\par\noindent
eiichi.nakai.math@vc.ibaraki.ac.jp
\par\noindent

\end{document}